# A New Theorem on N-th Order Differential Equation with Retarded Argument

Erdoğan Şen

*Namik Kemal University, Mathematics Department, 59030, Tekirdağ, Turkey*
*erdogan.maths@yahoo.com*

**Abstract.** In this work, using the well-known mean-value theorem (Lagrange's theorem) we obtain an inequality for n-th order differential equations with retarded argument. If the retarded argument vanishes then the inequality turns to an inequality for n-th order ordinary differential equations which plays a vital role in the proof of existence theorem for n-th order ordinary differential equations.

## INTRODUCTION

N-th order differential equations with deviating arguments are equations of the form

$$F\left(t, x(t), ..., x^{(m_0)}(t), x(t-\Delta_1(t)), ..., x^{(m_1)}(t-\Delta_1(t)), ..., x(t-\Delta_n(t)), ..., x^{(m_n)}(t-\Delta_n(t))\right) = 0$$

where $\Delta_i(t) \geq 0$, $i = 1, ..., n$, and $\max_{0 \leq i \leq n} m_i = n$. By $x^{(m_i)}(t-\Delta_i(t))$ we mean the $m_i$ th derivative of the function $x(z)$ evaluated at the point $z = t - \Delta_i(t)$.

A natural classification of equations with deviating arguments was proposed by G.A. Kamenskii [1]. When equation is solved for $x^{(m_0)}(t)$, it becomes

$$x^{(m_0)}(t) = f\left(t, x(t), x^{(m_0-1)}(t), x(t-\Delta_1(t)), ..., x^{(m_1)}(t-\Delta_1(t)), ..., x(t-\Delta_n(t)), ..., x^{(m_n)}(t-\Delta_n(t))\right).$$

Let $\eta = m_0 - \mu$, where $\mu = \max_{1 \leq i \leq n} m_i$. Equations for which $\eta > 0$ are called equations with retarded argument, those for which $\eta = 0$ are called equations of neutral type and those for which $\eta < 0$ are called equations of advanced type.

In this study we consider the n-th order differential equations with retarded argument (DERA)

$$L(y) = y^{(n)}(t) + M_1(t) y^{(n-1)}(t-\Delta(t)) + ... + M_n(t) y(t-\Delta(t)) = 0 \tag{1}$$

on an interval $I \subset \mathbb{R}$. Here $M_j(t)$, $\Delta(t)$ are continuous functions on $I$ for $j = 1, 2, ..., n$ and $1 \geq \Delta(t) \geq 0$ for each $t \in I$. In this work, by virtue of the well-known mean-value theorem (Lagrange's theorem) we obtain an inequality for the Equation (1). Similar inequality for second order DERA can be found in [2].

# AN INEQUALITY FOR N-TH ORDER DERA

**Theorem.** *Let us denote by every point with $t_{k_i}$ which is satisfying the mean-value theorem for a continuous solution $w^{(i)}(t_k)$ of (1) on $[t_k - \Delta(t_k), t_k] \subseteq I$ for each $t_k \in I$, $k \in J$ and $i = 0,1,...,n-1$ where $J$ is an index set. Also let us assume that $\sup_{t \in I} M_j(t) = M_{0_j}$, $j = 1,2,...,n$ where $M_{0_j}$s are real numbers. Then for all $t_k$ in $I$*

$$\left\| w(t_{k_i}) \right\| e^{-\psi \left| t_k - t_{k_i} \right|} \leq \left\| w(t_k) \right\| \leq \left\| w(t_{k_i}) \right\| e^{\psi \left| t_k - t_{k_i} \right|} \tag{2}$$

*where*

$$\left\| w(t_k) \right\| = \left[ \left| w(t_k) \right|^2 + \left| w'(t_k) \right|^2 + ... + \left| w^{(n-1)}(t_k) \right|^2 \right]^{1/2},$$

$$\psi = 1 + \left| M_{0_1} \right| (1 + \left| w^{(n)}(t_{k_{n-1}}) \right|) + \left| M_{0_2} \right| (1 + \left| w^{(n-1)}(t_{k_{n-2}}) \right|) + ... + \left| M_{0_n} \right| (1 + \left| w'(t_{k_0}) \right|).$$

**Proof.** In the case $\Delta(t) = 0$ the theorem can be proved as being in the theory of ordinary differential equations [3]. Now let us consider the case $\Delta(t) > 0$. From the mean-value theorem we can write the followings:

$$\frac{w(t_k) - w(t_k - \Delta(t_k))}{\Delta(t_k)} = w'(t_{k_0}),$$

$$\frac{w'(t_k) - w'(t_k - \Delta(t_k))}{\Delta(t_k)} = w''(t_{k_1}),$$

$$\vdots$$

$$\frac{w^{(n-1)}(t_k) - w^{(n-1)}(t_k - \Delta(t_k))}{\Delta(t_k)} = w^{(n)}(t_{k_{n-1}}).$$

Here $t_{k_i} \in \left( t_k - \Delta(t_k), t_k \right)$, $i = 0,1,...,n-1$. Thus

$$w\left( t_k - \Delta(t_k) \right) = w(t_k) - w'\left( t_{k_0} \right) \Delta(t_k),$$

$$w'\left( t_k - \Delta(t_k) \right) = w'(t_k) - w''\left( t_{k_1} \right) \Delta(t_k),$$

$$\vdots$$

$$w^{(n-1)}\left( t_k - \Delta(t_k) \right) = w^{(n-1)}(t_k) - w^{(n)}\left( t_{k_{n-1}} \right) \Delta(t_k),$$

and

$$|w(t_k - \Delta(t_k))| \leq |w(t_k)| + |w'(t_{k_0})|,$$
$$|w'(t_k - \Delta(t_k))| \leq |w'(t_k)| + |w''(t_{k_1})|,$$
$$\vdots \qquad (3)$$
$$|w^{(n-1)}(t_k - \Delta(t))| \leq |w^{(n-1)}(t_k)| + |w^{(n)}(t_{k_{n-1}})|.$$

Now we let $u(t_k) = \|w(t_k)\|^2$. Thus

$$u = w\overline{w} + w'\overline{w'} + \ldots + w^{(n-1)}\overline{w^{(n-1)}},$$

where $\overline{w(t_k)} = \overline{w(t_k)}$. Then

$$u' = w'\overline{w} + w''\overline{w'} + \ldots + w^{(n)}\overline{w^{(n-1)}} + w\overline{w'} + w'\overline{w''} + \ldots + w^{(n-1)}\overline{w^{(n)}}.$$

From the definition of a derivative it follows that $\overline{w'} = \overline{w}'$. Also $|w(t_k)| = |\overline{w(t_k)}|$. Therefore

$$|u'(t_k)| \leq 2|w(t_k)||w'(t_k)| + 2|w'(t_k)||w''(t_k)| + \ldots + 2|w^{(n-1)}(t_k)||w^{(n)}(t_k)|. \qquad (4)$$

Since $w$ satisfies $L(w) = 0$ we have

$$w^{(n)}(t_k) = -\left[M_1(t_k)w^{(n-1)}(t_k - \Delta(t_k)) + M_2(t_k)w^{(n-2)}(t_k - \Delta(t_k)) + \ldots + M_n(t_k)w(t_k - \Delta(t_k))\right]$$

and hence applying (3)

$$|w^{(n)}(t_k)| \leq |M_{0_1}|\left\{|w^{(n-1)}(t_k)| + |w^{(n)}(t_{k_{n-1}})|\right\} + |M_{0_2}|\left\{|w^{(n-2)}(t_k)| + |w^{(n-1)}(t_{k_{n-2}})|\right\} + \ldots + |M_{0_n}|\left\{|w(t_k)| + |w'(t_{k_0})|\right\}. \qquad (5)$$

Using (5) in the equation (4) we obtain

$$|u'(t_k)| \leq 2|w(t_k)||w'(t_k)| + 2|w'(t_k)||w''(t_k)| + \ldots + 2|w^{(n-2)}(t_k)||w^{(n-1)}(t_k)| + 2|M_{0_1}||w^{(n-1)}(t_k)|^2$$
$$+ 2|M_{0_2}||w^{(n-2)}(t_k)||w^{(n-1)}(t_k)| + \ldots + 2|M_{0_{n-1}}||w'(t_k)||w^{(n-1)}(t_k)|$$
$$+ 2|M_{0_n}||w(t_k)||w^{(n-1)}(t_k)| + 2|M_{0_1}||w^{(n)}(t_{k_{n-1}})||w^{(n-1)}(t_k)|$$
$$+ 2|M_{0_2}||w^{(n-1)}(t_{k_{n-2}})||w^{(n-1)}(t_k)| + \ldots + 2|M_{0_n}||w'(t_{k_0})||w^{(n-1)}(t_k)|.$$

Now applying the elementary fact that if $f$ and $g$ are any two functions then

$$2|f||g| \leq |f|^2 + |g|^2,$$

we get

$$|u'(t_j)| \leq (1 + |M_{0_n}|)|w(t_k)|^2 + (2 + |M_{0_{n-1}}|)|w'(t_k)|^2 + \ldots + (2 + |M_{0_2}|)|w^{(n-2)}(t_k)|^2 +$$
$$(1 + 2|M_{0_1}| + |M_{0_2}| + \ldots + |M_{0_n}| + 2[|M_{0_1}||w^{(n)}(t_{k_{n-1}})| + |M_{0_2}||w^{(n-1)}(t_{k_{n-2}})| + \ldots + |M_{0_n}||w'(t_0)|])|w^{(n-1)}(t_k)|^2$$

or

$$|u'(t_k)| \leq 2\psi u(t_k).$$

This is equivalent to

$$-2\psi u(t_k) \leq u'(t_k) \leq 2\psi u(t_k). \tag{6}$$

And these inequalities lead directly to (2). Indeed consider the right inequality which can be written as $u' - 2\psi u \leq 0$. It is equivalent to

$$e^{-2\psi t_k}(u' - 2\psi u) = (e^{-2\psi t_k} u)' \leq 0.$$

If $t_k > t_{k_i}$ we integrate from $t_k$ to $t_{k_i}$ obtaining $e^{-2\psi t_k} u(t_k) - e^{-2\psi t_{k_i}} u(t_{k_i}) \leq 0$ or

$$u(t_k) \leq u(t_{k_i}) e^{2\psi(t_k - t_{k_i})}. \tag{7}$$

If we take square root of each side in (7), we get

$$\|w(t_k)\| \leq \|w(t_{k_i})\| e^{\psi(t_k - t_{k_i})}, \quad t_k > t_{k_i}.$$

The left inequality in (6) similarly implies

$$\|w(t_{k_i})\| e^{-\psi(t_k - t_{k_i})} \leq \|w(t_k)\|, \quad t_k > t_{k_i}.$$

And therefore

$$\|w(t_{k_i})\| e^{-\psi|t_k - t_{k_i}|} \leq \|w(t_k)\| \leq \|w(t_{k_i})\| e^{\psi|t_k - t_{k_i}|}, \quad t_k > t_{k_i}$$

which is just (2) for $t_k > t_{k_i}$. The case $t_k < t_{k_i}$ may be considered analogically. Thus the theorem is proved.

We also want to note that in the case of $\Delta(t) = 0$ inequality (2) turns to an inequality for n-th order ordinary differential equations which plays a vital role in the proof of existence theorem for n-th order ordinary differential equations (see [3]).